\documentclass{article}
\pdfoutput=1

\usepackage{arxiv}
\usepackage{natbib}

\usepackage[utf8]{inputenc} 
\usepackage[T1]{fontenc}    
\usepackage{hyperref}       
\usepackage{url}            
\usepackage{booktabs}       
\usepackage{amsfonts}       
\usepackage{nicefrac}       
\usepackage{microtype}      
\usepackage{xcolor}         

\usepackage{mathtools}
\usepackage{nicematrix}
\usepackage{makecell}
\usepackage{amsthm,amssymb}
\usepackage{algorithm}
\usepackage{algpseudocode}
\usepackage{cleveref}
\usepackage{multirow}
\usepackage{enumitem}
\allowdisplaybreaks[3]

\newtheorem{theorem}{Theorem}
\newtheorem{lemma}{Lemma}

\newtheorem{assumption}{Assumption}

\newtheorem{definition}{Definition}
\crefname{assumption}{assumption}{assumptions}

\newcommand{\R}{\mathbb{R}}

\def\<#1,#2>{\langle #1,#2\rangle}

\newcommand{\norm}[1]{\|#1\|}
\newcommand{\sqn}[1]{\norm{#1}^2}

\newcommand{\cO}{\mathcal{O}}

\DeclareMathOperator*{\argmin}{arg\,min}
\DeclareMathOperator*{\Argmin}{Arg\,min}

\title{The First Optimal Acceleration of High-Order Methods in Smooth Convex Optimization}

\author{%
	Dmitry~Kovalev\\
	KAUST\thanks{King Abdullah University of Science and Technology, Thuwal, Saudi Arabia}\\
	\texttt{dakovalev1@gmail.com} 
	\And
	Alexander~Gasnikov\\
	IITP~RAS\thanks{Institute for Information Transmission Problems RAS, Moscow, Russia}\\
	\texttt{gasnikov@yandex.ru}
}

\begin{document}

\maketitle

\begin{abstract}

In this paper, we study the fundamental open question of finding the optimal high-order algorithm for solving smooth convex minimization problems. \citet{arjevani2019oracle} established the lower bound $\Omega\left(\epsilon^{-2/(3p+1)}\right)$ on the number of the $p$-th order oracle calls required by an algorithm to find an $\epsilon$-accurate solution to the problem, where the $p$-th order oracle stands for the computation of the objective function value and the derivatives up to the order $p$. However, the existing state-of-the-art high-order methods of \citet{gasnikov2019near,bubeck2019near,jiang2019optimal} achieve the oracle complexity $\mathcal{O}\left(\epsilon^{-2/(3p+1)} \log (1/\epsilon)\right)$, which does not match the lower bound. The reason for this is that these algorithms require performing a complex binary search procedure, which makes them neither optimal nor practical. We fix this fundamental issue by providing the first algorithm with $\mathcal{O}\left(\epsilon^{-2/(3p+1)}\right)$ $p$-th order oracle complexity.

\end{abstract}

\section{Introduction}

Let $\R^d$ be a finite-dimensional Euclidean space and let $f(x)\colon \R^d \rightarrow \R$ be a convex, $p$ times continuously differentiable function with $L_p$-Lipschitz $p$-th order derivatives. Our goal is to solve the following convex minimization problem:
\begin{equation}\label{eq:main}
	f^* = \min_{x \in \R^d} f(x).
\end{equation}
In this work, we assume access to the $p$-th order oracle associated  with function $f(x)$. That is, given an arbitrary point $x \in \R^d$, we can compute the function value and the derivatives of function $f(x)$ up to order $p$.

\paragraph{First-order methods.}
When $p=1$, first-order methods, such as gradient descent, are typically used for solving problem~\eqref{eq:main}.
The lower bound $\Omega(\epsilon^{-1/2})$ on the number of the gradient evaluations required by these algorithms to find an $\epsilon$-accurate solution was established by \citet{nemirovskij1983problem,nesterov2003introductory}, while the optimal algorithm matching this lower bound is called Accelerated Gradient Descent and was developed by \citet{nesterov1983method}.


\paragraph{Second-order methods.}
In contrast to the first-order methods, the understanding of the second-order methods ($p=2$) was developed relatively recently. \citet{nesterov2006cubic} developed the cubic regularized variant of Newton's method. This algorithm achieves global convergence with the oracle complexity $\cO(\epsilon^{-1/2})$, which cannot be achieved with the standard Newton's method. \citet{nesterov2008accelerating} also developed an accelerated version of the cubic regularized Newton's method with $\cO\left( \epsilon^{-1/3}\right)$ second-order oracle complexity. A few years later, \citet{monteiro2013accelerated} developed the Accelerated Hybrid Proximal Extragradient (A-HPE) framework and combined it with a trust region Newton-type method. The resulting algorithm, called Accelerated Newton Proximal Extragradient (A-NPE), achieved the second-order oracle complexity of $\cO\left(\epsilon^{-2/7}\log (1/\epsilon) \right)$. In 2018, \citet{arjevani2019oracle} established the lower bound $\Omega \left(\epsilon^{-2/7}\right)$ on the number of the second-order oracle calls required by an algorithm to find an $\epsilon$-accurate solution\footnote{There is also a work of \citet{agarwal2018lower}, which provides the lower complexity bounds for high-order optimization. However, their lower bounds are worse than the lower bounds of \citet{arjevani2019oracle}.}, which is almost achieved by the A-NPE algorithm of \citet{monteiro2013accelerated}, up to the logarithmic factor $\log (1/\epsilon)$. However, the optimal second-order algorithms for solving smooth convex minimization problems remain to be unknown.

\paragraph{High-order methods.}
In the case when $p>2$, the situation is very similar to the second-order case. \citet{nesterov2021implementable} developed the generalization of the cubic regularized Newton method to the high-order case and called them tensor methods. \citet{nesterov2021implementable} provided both non-accelerated and accelerated $p$-th order tensor methods with the oracle complexity $\cO\left(\epsilon^{-1/p}\right)$ and $\cO\left(\epsilon^{-1/(p+1)}\right)$, respectively.\footnote{\citet{nesterov2021implementable} also provided the lower complexity bounds that coincide with the lower bounds of \citet{arjevani2019oracle}.} Later, three independent groups of researchers \citep{gasnikov2019optimal,bubeck2019near,jiang2019optimal} used the A-HPE framework to develop the near-optimal tensor methods with the oracle complexity $\cO\left(\epsilon^{-2/(3p+1)}\log (1/\epsilon)\right)$. Similarly to the case $p=2$, these algorithms match the lower complexity bound $\Omega\left(\epsilon^{-2/(3p+1)}\right)$ of \citet{arjevani2019oracle}, up to the logarithmic factor $\log (1/\epsilon)$.

\subsection{Main Contribution: Optimal Second-Order and High Order Methods}

\begin{table}[t]
	\caption{Comparison of the first-order, second-order and high-order methods for smooth convex optimization in the oracle complexities (see \Cref{def:complexity}), which depend on the smoothness constant $L_p$ (see \Cref{ass:lip}), the distance to the solution $R$ (see \Cref{ass:dist}), and the accuracy $\epsilon$ (see \Cref{def:accuracy}).}
	\label{tab:main}
	\centering
{
		\begin{NiceTabular}{|c|c|Wc{6em}|}
			\toprule
			\bf Algorithm Reference & \bf Oracle Complexity& \bf  Order  \\
			\midrule
			\makecell{\citet{nesterov1983method}} & $\cO\left(\left({L_1R^2}/{\epsilon}\right)^{1/2}\right)$&
			\Block{2-1}{\bf\makecell{First-Order\\Methods\\($p=1$)}}\\
			\cmidrule{1-2}
			\makecell{Lower Bound \citep{nemirovskij1983problem}} & $\Omega\left(\left({L_1R^2}/{\epsilon}\right)^{1/2}\right)$&\\
			\midrule
			\makecell{\citet{nesterov2006cubic}}&$\cO\left(\left({L_2R^3}/{\epsilon}\right)^{1/2}\right)$ & \Block{5-1}{\bf\makecell{Second-Order\\Methods\\($p=2$)}}\\
			\cmidrule{1-2}
			\makecell{ \citet{nesterov2008accelerating}} & $\cO\left(\left({L_2R^3}/{\epsilon}\right)^{1/3}\right)$&\\
			\cmidrule{1-2}
			\makecell{\citet{monteiro2013accelerated}} & $\cO\left(\left({L_2R^3}/{\epsilon}\right)^{2/7} \log ({1}/{\epsilon})\right)$&\\
			\cmidrule{1-2}
			\makecell{\bf \Cref{alg:ot} \bf (This Paper)}&$\cO\left(\left({L_2R^3}/{\epsilon}\right)^{2/7} \right)$&\\
			\cmidrule{1-2}
			\makecell{Lower Bound \citep{arjevani2019oracle}} & $\Omega\left(\left({L_2R^3}/{\epsilon}\right)^{2/7} \right)$\\
			\midrule
			\makecell{ \citet{nesterov2021implementable}} & $\cO\left(\left({L_pR^{p+1}}/{\epsilon}\right)^{1/p}\right)$& \Block{5-1}{\bf\makecell{High-Order\\Methods\\($p\geq 2$)}} \\
			\cmidrule{1-2}
			\makecell{\citet{nesterov2021implementable}} & $\cO\left(\left({L_pR^{p+1}}/{\epsilon}\right)^{1/{p+1}}\right)$&\\
			\cmidrule{1-2}
			\makecell{\citet{gasnikov2019near}} & $\cO\left(\left({L_pR^{p+1}}/{\epsilon}\right)^{2/{(3p+1)}} \log ({1}/{\epsilon})\right)$&\\
			\cmidrule{1-2}
			\makecell{\bf \Cref{alg:ot} \bf (This Paper)}&$\cO\left(\left({L_pR^{p+1}}/{\epsilon}\right)^{2/{(3p+1)}} \right)$\\
			\cmidrule{1-2}
			\makecell{Lower Bound  \citep{arjevani2019oracle}} & $\Omega\left(\left({L_pR^{p+1}}/{\epsilon}\right)^{2/(3p+1)} \right)$&\\
			\bottomrule
	\end{NiceTabular}}
\end{table}

The review of the second-order and high-order methods that we made above identifies the following fundamental open question:
\begin{center}
	\em Can we design an optimal $p$-th order algorithm ($p\geq 2$) for solving smooth convex minimization problems with the oracle complexity matching the lower bounds?
\end{center}
The lack of an answer to this question reveals a significant gap in the understanding of the high-order optimization compared to the first-order optimization.
We give a positive answer to this question. That is, we provide the first optimal high-order algorithm with the $p$-th order oracle complexity $\cO\left(\left(\epsilon^{-2/(3p+1)}\right)\right)$ that matches the lower bounds of \citet{arjevani2019oracle}. This is the main contribution of our work.

Our paper is organized as follows: {\bf(i)}~in \Cref{sec:prelim}, we briefly introduce the tensor approximations and provide necessary assumptions and definitions; {\bf(ii)}~in \Cref{sec:not}, we describe the existing near-optimal high-order methods and identify their main flaws that prevent them from being optimal and practical algorithms; {\bf(iii)}~in \Cref{sec:ot}, we describe the development of our optimal high-order algorithm and provide its theoretical convergence analysis.



\section{Preliminaries}\label{sec:prelim}

By $\norm{\cdot} \colon \R^d \rightarrow \R$ and $\<\cdot,\cdot> \colon \R^d \times \R^d \rightarrow \R$, we denote the standard Euclidean norm and scalar product on $\R^d$.
Given a $p$ times continuously differentiable function $g(x) \colon R^d \rightarrow \R$ and index $i \in \{1,2,\ldots, p\}$, by $\nabla^i g(x)[h]^i \colon \R^d \rightarrow \R$ we denote the following homogeneous polynomial:
\begin{equation}
	\nabla^i g(x) [h]^i = \sum_{j_1,\ldots,j_i = 1}^d \frac{\partial^i g}{\partial x_{j_1}\cdots \partial x_{j_i} }(x) \cdot h_{j_1}\cdots h_{j_i},
\end{equation}
where $x = (x_1,\ldots,x_d)\in \R^d$, $h = (h_1,\ldots,h_d) \in \R^d$, and 
\begin{equation}
	\frac{\partial^i g}{\partial x_{j_1}\cdots \partial x_{j_i} }(x)
\end{equation}
is the $i$-th order partial derivative of function $g(x)$ at point $x$ with respect to variables $x_{j_1},\ldots,x_{j_i}$.
 For instance, if $i=1$, then $\nabla^1 g(x)[h] = \<\nabla g(x), h>$, where $\nabla g(x) \in \R^d$ is the gradient of function $g(x)$; if $i=2$, then $\nabla^2 g(x)[h] =\< \nabla^2 f(x) h,h>$, where $\nabla^2 f(x) \in \R^{d\times d}$ is the Hessian of function $f(x)$.
We can write the $p$-th order Taylor approximation of function $g(x)$ at point $z \in \R^d$:
\begin{equation}
	\Phi_g^p(x;z) = g(z) + \sum_{i=1}^p \frac{1}{i!} \nabla^i g(z)[x-z],
\end{equation}
It is well known that the Taylor polynomial $\Phi_g^p(x;z)$ approximates function $g(x)$, if point $x$ is close enough to point $z$:
\begin{equation}
	g(x) = \Phi_g(x;z) + R_g^p(x;z)\norm{x-z}^p,
\end{equation}
where $R_g^p(\cdot;z)\colon \R^d\rightarrow \R$ is a function that satisfies $\lim_{x \rightarrow z} R_g^p(x;z) = 0$.

As mentioned earlier, we assume that the objective function $f(x)$ of the main problem~\eqref{eq:main} is $p$ times continuously differentiable and has $L_p$-Lipschitz $p$-th order derivatives. It is formalized via the following definition.
\begin{assumption}\label{ass:lip}
	Function $f(x)$ is $p$-times continuously differentiable, convex, and has $L_p$-Lipschitz $p$-th order derivatives, i.e., for all $x_1,x_2 \in \R^d$ the following inequality holds:
	\begin{equation*}
		\max \{ \lvert\nabla^p f(x_1)[h] - \nabla^p f(x_2)[h] \rvert :  h \in \R^d, \norm{h} \leq 1 \} \leq L_p \norm{x_1 - x_2}.
	\end{equation*}
\end{assumption}
Theorem~1 of \citet{nesterov2021implementable} implies that under \Cref{ass:lip}, function $f(x)$ has the following convex upper bound:
\begin{equation}\label{eq:upper}
	f(x) \leq \Phi_g^p(x;z) + \frac{pM}{(p+1)!}\norm{x-z}^{p+1},
\end{equation}
where $M \geq L_p$ and $z \in \R^d$. Hence, an obvious approach to solving problem~\eqref{eq:main} is to perform the minimization of this upper bound instead of minimizing the function $f(x)$. This approach naturally leads to the following iterative process:
\begin{equation}\label{eq:tensor}
	x^{k+1} \in \Argmin_{x \in \R^d} \Phi_g^p(x;x^k) + \frac{pM}{(p+1)!}\norm{x-x^k}^{p+1}.
\end{equation}
In the case $p=2$, this iterative process is known as the cubic regularized Newton's method of \citet{nesterov2006cubic}, and in the case $p>2$, it is known as the tensor method of \citet{nesterov2021implementable}. Minimization procedures similar to~\eqref{eq:tensor} are widely used in high-order optimization methods. It will also be used in the development of our optimal algorithm.

We also have the following assumption which requires problem~\eqref{eq:main} to have at least a single solution $x^* \in \R^d$. It is a standard assumption for the majority of works on convex optimization. 
\begin{assumption}\label{ass:dist}
	There exists a constant $R > 0$ and at least a single solution $x^*$ to problem~\eqref{eq:main}, such that $\norm{x^0 - x^*} \leq R$, where $x^0 \in \R^d$ is the starting point that we use as an input for a given algorithm for solving the problem.
\end{assumption}

Finally, we have the following definitions that formalize the notions of $\epsilon$-accurate solution of a problem, $p$-th order oracle call, and oracle complexity of an algorithm.
\begin{definition}\label{def:accuracy}
	We call vector $\hat x\in \R^d$ an $\epsilon$-accurate solution of problem~\eqref{eq:main}, if for a given accuracy $\epsilon > 0$ it satisfies $f(\hat x) - f^* \leq \epsilon$.
\end{definition}
\begin{definition}
	Given an arbitrary vector $x \in \R^d$ by the $p$-th order oracle call at $x$, we denote the computation of the function value $f(x)$ and the derivatives $\nabla^1 f(x)[\cdot], \ldots, \nabla^p f(x)[\cdot]$.
\end{definition}
\begin{definition}\label{def:complexity}
	By the $p$-th order oracle complexity of a $p$-th order algorithm for solving problem~\eqref{eq:main}, we denote the number of $p$-th order oracle calls required by the algorithm to find an $\epsilon$-accurate solution of the problem for a given accuracy $\epsilon > 0$.
\end{definition}

\section{Near-Optimal Tensor Methods}\label{sec:not}

In this section, we revisit the state-of-the-art high-order optimization algorithms that include the A-NPE method of \citet{monteiro2013accelerated} in the $p=2$ case and the near-optimal tensor methods of \citet{gasnikov2019optimal,bubeck2019near,jiang2019optimal} in the general $p>2$ case. We start with describing the key ideas behind the development of these algorithms to understand how they work. Then, we identify the main flaws of the algorithms that prevent them from being optimal and practical.

Note that the A-NPE method and near-optimal tensor methods have the following substantial similarities:
{\bf(i)}~both the A-NPE and near-optimal tensor methods are based on the A-HPE framework of \citet{monteiro2013accelerated};
{\bf(ii)}~the oracle complexity of the near-optimal tensor methods recovers the oracle complexity of the A-NPE method in the case $p=2$; {\bf (iii)} these algorithms have the same issue: the requirement to perform the complex binary search procedure at each iteration which makes them neither optimal nor practical.
Hence, we will further leave out the description of the A-NPE method of \citet{monteiro2013accelerated} and consider only the near-optimal tensor methods of \citet{gasnikov2019optimal,bubeck2019near,jiang2019optimal}.


\subsection{A-HPE Framework}

\begin{algorithm}[t]
	\caption{A-HPE Framework}
	\label{alg:hpe}
	\begin{algorithmic}[1]
		\State {\bf input: } $x^0 = x_f^0 \in \R^d$
		\State {\bf parameters:} $\sigma \in [0,1]$, $K \in \{1,2,\ldots\}$
		\State $\beta_{-1} = 0$
		\For{ $k = 0,1,2,\ldots, K-1$}
		\State compute $x_f^{k+1}\in \R^d$, $\lambda_k > 0$ such that \label{hpe:line:aux}
		\begin{equation}\label{hpe:aux}
						\norm{\nabla f(x_f^{k+1}) + \lambda_{k}^{-1}(x_f^{k+1} - x_g^k)} \leq \sigma \lambda_k^{-1}\norm{x_f^{k+1} - x_g^k},
		\end{equation}
		\hskip\algorithmicindent where $x_g^k \in \R^d$ and $\alpha_k \in (0,1]$ are defined as
		\begin{equation}\label{hpe:aux2}
			x_g^k = \alpha_k x^k + (1-\alpha_k)x_f^k, \quad \alpha_k = \eta_k/\beta_k, 
		\end{equation}
		\hskip\algorithmicindent and $\eta_k > 0$ and $\beta_k > 0$ are defined by the following system:
		\begin{equation}\label{hpe:aux3}
			\beta_{k-1} + \eta_k = \beta_k, \quad \beta_k \lambda_k = \eta_k^2. 
		\end{equation} 
		\State $x^{k+1} = x^k - \eta_k \nabla f(x_f^{k+1})$ \label{hpe:line:grad}
		\EndFor
		\State {\bf output:} $x_f^{K}$
	\end{algorithmic}
\end{algorithm}

The main component in the development of the near-optimal tensor methods of \citet{gasnikov2019optimal,bubeck2019near,jiang2019optimal} is the Accelerated Hybrid Proximal Extragradient (A-HPE) framework of \citet{monteiro2013accelerated}. This algorithmic framework can be seen as a generalization of the Accelerated Gradient Descent of \citet{nesterov1983method}. It is formalized as \Cref{alg:hpe}.
Next, we recall the main theorem by \citet{monteiro2013accelerated}, which describes the convergence properties of \Cref{alg:hpe}.

\begin{theorem}[\citet{monteiro2013accelerated}]\label{hpe:thm}
	The iterations of \Cref{alg:hpe} satisfy the following inequality:
	\begin{equation}\label{hpe:thm:eq}
		2\beta_{K-1} (f(x_f^K) - f^*)  + (1-\sigma^2)\sum_{k=0}^{K-1} \alpha_{k}^{-2}\sqn{x_f^{k+1} - x_g^k} \leq R^2.
	\end{equation}
\end{theorem}

Note that \Cref{alg:hpe} requires finding $x_f^{k+1}$ satisfying condition~\eqref{hpe:aux} on \cref{hpe:line:aux}. This condition can be rewritten as follows:
\begin{equation}\label{hpe:monteiro-svaiter}
	\norm{\nabla A_{\lambda_k} (x_f^{k+1};x_g^k)} \leq \sigma\lambda_{k}^{-1}\norm{x_f^{k+1} - x_g^k},
\end{equation}
where function $A_{\lambda}(\cdot;z) \colon \R^d \rightarrow \R$ for $\lambda > 0$ and $z \in \R^d$ is defined as
\begin{equation}\label{eq:A}
	A_{\lambda}(x;z) = f(x) + \frac{1}{2\lambda}\sqn{x - z}.
\end{equation}

\subsection{Application to High-Order Minimization}

In order to perform the computation on \cref{hpe:line:aux} of \Cref{alg:hpe}, we need to find $x_f^{k+1} \in \R^d$ that satisfies condition~\eqref{hpe:aux}. As we mentioned earlier, condition~\eqref{hpe:aux} is equivalent to \eqref{hpe:monteiro-svaiter}, which involves the gradient norm $\norm{\nabla A_{\lambda_k}(\cdot;x_g^k)}$ at point $x_f^{k+1}$. Function $A_{\lambda_k}(\cdot;x_g^k)$ has $L_p$-Lipschitz $p$-th order derivatives for $p \geq 2$ due to its definition~\eqref{eq:A} and \Cref{ass:lip}.\footnote{$\nabla^p A(x;x_g^k)[h] = \nabla^p f(x)[h]$ when $p>2$, and $\nabla^2 A(x;x_g^k)[h] = \nabla^2 f(x)[h] + \lambda^{-1}\sqn{h}$.}
Hence, it has the following upper bound, thanks to Theorem~1 of \citet{nesterov2021implementable}:
\begin{equation}\label{eq:upper2}
	A_{\lambda_k} (x;x_g^k) \leq \Phi_{A_{\lambda_k} (\cdot;x_g^k)}^p (x;x_g^k) + \frac{pM}{(p+1)!}\norm{x-x_g^k}^{p+1}.
\end{equation}
It turns out that $x_f^{k+1}$ can be obtained by minimizing this upper bound:
\begin{equation}\label{eq:tensor_step}
	x_f^{k+1} = \argmin_{x \in \R^d}\Phi_{A_{\lambda_k} (\cdot;x_g^k)}^p (x;x_g^k) + \frac{pM}{(p+1)!}\norm{x-x_g^k}^{p+1},
\end{equation}
where $M > L_p$.\footnote{We require the strict inequality to ensure the uniform convexity of upper bound~\eqref{eq:upper2}, which implies the uniqueness and the existence of the minimizer in~\eqref{eq:tensor_step}.}
Indeed, by Lemma~1 of \citet{nesterov2021implementable}, we have
\begin{equation}
	\norm{\nabla A_{\lambda_k}(x_f^{k+1};x_g^k)} \leq \frac{pM + L_p}{p!}\norm{x_f^{k+1} - x_g^k}^{p}.
\end{equation}
Hence, to satisfy condition~\eqref{hpe:monteiro-svaiter}, we choose $\lambda_k$ in the following way:
\begin{equation}\label{eq:search}
	\frac{\sigma p!}{2(pM + L_p)}\norm{x_f^{k+1} - x_g^k}^{1-p} \leq \lambda_k \leq \frac{\sigma p!}{(pM + L_p)}\norm{x_f^{k+1} - x_g^k}^{1-p}.
\end{equation}
Here, the upper bound on $\lambda_k$ ensures condition~\eqref{hpe:monteiro-svaiter}, while the lower bound prevents stepsize $\lambda_k$ from being too small, which would hurt the convergence rate.
The resulting near-optimal tensor method is formalized as \Cref{alg:not}.
It has the following convergence rate:
\begin{equation}
	f(x_f^K) - f^* \leq  \frac{\mathrm{const} \cdot L_p\norm{x^0 - x^*}^{p+1}}{K^{\frac{3p+1}{2}}},
\end{equation}
where $K$ is the number of iterations.
The proof of this convergence rate involves condition~\eqref{eq:search} and \Cref{hpe:thm}. It is given in the works of \citet{gasnikov2019optimal,bubeck2019near,jiang2019optimal}.

\begin{algorithm}[t]
	\caption{Near-Optimal Tensor Method}
	\label{alg:not}
	\begin{algorithmic}[1]
		\State {\bf input: } $x^0 = x_f^0 \in \R^d$
		\State {\bf parameters:} $M > 0$, $K \in \{1,2,\ldots\}$
		\State $\beta_{-1} = 0$
		\For{ $k = 0,1,2,\ldots, K-1$}
		\State compute
		$
		\begin{cases}		
			\lambda_k > 0& \text{satisfying } \eqref{eq:search}\\
			x_f^{k+1} \in \R^d & \text{satisfying } \eqref{eq:tensor_step}\\
			x_g^k \in \R^d,\alpha_k \in (0,1]& \text{satisfying } \eqref{hpe:aux2}\\
			\eta_k,\beta_k > 0& \text{satisfying } \eqref{hpe:aux3}
		\end{cases}
		$ \label{not:line:aux}
		\State $x^{k+1} = x^k - \eta_k \nabla f(x_f^{k+1})$
		\EndFor
		\State {\bf output:} $x_f^{K}$
	\end{algorithmic}
\end{algorithm}

\subsection{The Problems with the Existing Algorithms}

\Cref{alg:not} requires finding $\lambda_k$ satisfying condition~\eqref{eq:search} at each iteration. According to \cref{not:line:aux} of \Cref{alg:not}, $\lambda_k$ depends on $x_f^{k+1}$ via~\eqref{eq:search}, which depends on $x_g^k$ via~\eqref{eq:tensor_step}, which depends on $\eta_k,\beta_k$ via~\eqref{hpe:aux2}, which depend on $\lambda_k$ via~\eqref{hpe:aux3}. Hence, computation of stepsize $\lambda_k$ depends on $\lambda_k$ itself and there is no explicit way to perform the computation on \cref{not:line:aux}. 

The algorithms of \citet{gasnikov2019optimal,bubeck2019near,jiang2019optimal} use various binary search procedures to find $\lambda_k$ and perform the computation on \cref{not:line:aux}. However, such procedures are costly and require many iterations to converge.
For instance, \citet{bubeck2019near} show that their variant of binary search requires the following number of $p$-th order oracle calls to find $\lambda_k$ satisfying condition~\eqref{eq:search}:
\begin{equation}
	\cO \left(\log \frac{L_pR^{p+1}}{\epsilon}\right).
\end{equation}
The same complexity (up to constant factors) for similar binary search procedures was established in the works of \citet{nesterov2021inexact,jiang2019optimal}, and in the work of \citet{monteiro2013accelerated} for the $p=2$ case. Hence, the total oracle complexity of \Cref{alg:not} is $\cO \left(  \epsilon^{-2/(3p+1)} \log (1/ \epsilon)\right)$
which does not match the lower bound of \citet{arjevani2019oracle}.

The additional logarithmic factor in the oracle complexity of \Cref{alg:not} raises the question whether it is superior to the accelerated tensor method of \citet{nesterov2021implementable} in practice. On the one hand, \citet{gasnikov2019optimal} provided an experimental study that showed the practical superiority of \Cref{alg:not} over the algorithm of \citet{nesterov2021implementable}. However, this experimental comparison is utterly unfair because it considers only the iteration complexity of the algorithms, which does not take into account the oracle complexity of the binary search procedure.

\section{The First Optimal Tensor Method}\label{sec:ot}

In the previous section, we described the main issues with the existing high-order methods that prevent them from being optimal and practical algorithms for solving problem~\eqref{eq:main}. In this section, we will show how to construct an algorithm that does not have those issues. More precisely, we will develop the first optimal $p$-th order algorithm ($p\geq 2$) for solving main problem~\eqref{eq:main}.

\subsection{The Key Idea}

The crucial mistake \citet{gasnikov2019optimal,bubeck2019near,jiang2019optimal} made while creating their algorithms is that they fixed the procedure of computing $x_f^{k+1}$ on \cref{hpe:line:aux} of \Cref{alg:hpe} using formula~\eqref{eq:tensor_step} and then developed the procedure for computing $\lambda_k$, which turned out to be inefficient. We will go the opposite way. That is, we choose parameters $\lambda_k$ in advance in such a way that they ensure the optimal convergence rate and then provide an efficient procedure for finding $x_f^{k+1}$ satisfying condition~\eqref{hpe:aux}.
Let $\eta_k$ be defined as follows:
\begin{equation}\label{eq:eta}
	\eta_k = \eta (1 + k)^{\frac{3p-1}{2}},
\end{equation}
where $\eta > 0$ is a parameter. Using~\eqref{hpe:aux3}, we can compute $\beta_k$ and $\lambda_k$ as follows:
\begin{equation}\label{eq:betalambda}
	\beta_k = \eta\sum_{l=0}^{k} (1 + l)^{\frac{3p-1}{2}}, \qquad \lambda^k = \frac{\eta (1+k)^{3p-1}}{\sum_{l=0}^{k} (1 + l)^{\frac{3p-1}{2}}}.
\end{equation}
The following lemma provides a lower bound on $\beta_k$ and an upper bound on $\lambda_k$.
\begin{lemma}\label{lem:betalambda}
	Parameters $\beta_k$ and $\lambda_k$ defined by~\eqref{eq:betalambda} satisfy the following inequalities:
	\begin{equation}
		\beta_k \geq \frac{2\eta}{(3p+1)} (k+1)^{\frac{3p+1}{2}}, \qquad 
		\lambda_k \leq \frac{\eta(3p+1)}{2}(1+k)^{\frac{3(p-1)}{2}}.
	\end{equation}
\end{lemma}
\Cref{lem:betalambda} and \Cref{hpe:thm} immediately imply the convergence rate $\cO(1/k^{(3p+1)/2})$, which matches the lower bound of \citet{arjevani2019oracle}. Hence, the only remaining question is how to compute $x_f^{k+1}$ satisfying~\eqref{hpe:aux} efficiently. To be precise, we need to develop a procedure that can perform this computation using $\cO(1)$ of $p$-th order oracle calls.

\subsection{Tensor Extragradient Method for Gradient Norm Reduction}

\begin{algorithm}[t]
	\caption{Tensor Extragradient Method}
	\label{alg:eg}
	\begin{algorithmic}[1]
		\State {\bf input:} $x^{k,0} = x_g^k \in \R^d$, $A^k(\cdot) =  A_{\lambda_{k}}(\cdot;x_g^k)$
		\State {\bf parameters: } $M > 0$
		\State $t=-1$
		\Repeat
		\State $t = t+1$
		\State compute $x^{k,t+1/2}\in\R^d$ as follows:
		\begin{equation}\label{eg:tensor_step}
			x^{k,t+1/2} = \argmin_{x \in \R^d} \Phi_{A^k}^p (x;x^{k,t}) + \frac{pM}{(p+1)!}\norm{x - x^{k,t}}^{p+1}
		\end{equation}  \label{eg:line:1}
		\State $x^{k,t+1} = x^{k,t} - \left(\frac{M\norm{x^{k,t+1/2} - x^{k,t}}^{p-1}}{(p-1)!}\right)^{-1} \nabla A^k(x^{k,t+1/2})$ \label{eg:line:2}
		\Until{$\norm{\nabla A^k(x^{k,t+1/2})} \leq \sigma \lambda_k^{-1}\norm{x^{k,t+1/2} - x^{k,0}}$}\label{eg:line:until}
		\State $T^k = t+1$
		\State {\bf output:} $x_f^{k+1} = x^{k,T^k-1/2}$
	\end{algorithmic}
\end{algorithm}

In this subsection, we develop an efficient procedure for computing $x_f^{k+1}$ satisfying condition~\eqref{hpe:aux}.
As we mentioned earlier, condition~\eqref{hpe:aux} is equivalent to~\eqref{hpe:monteiro-svaiter}, which is an upper bound on the gradient norm $\norm{\nabla A_{\lambda_k}(\cdot; x_g^k)}$ at point $x_f^{k+1}$. Hence, we need an algorithm for the gradient norm reduction in the following smooth high-order convex minimization problem:
\begin{equation}
	x^{k,*} = \argmin_{x \in \R^d} A_{\lambda_{k}}(x;x_g^k).
\end{equation}
In this subsection, we provide such an algorithm. We call the algorithm Tensor Extragradient Method. It is formalized as \Cref{alg:eg}. In the case $p=1$, this algorithm recovers the extragradient method of \citet{korpelevich1976extragradient}. \Cref{alg:eg} can be seen as a generalization of the extragradient method for high-order optimization.

One can observe that due to \cref{eg:line:until} of \Cref{alg:eg}, $x_f^{k+1} = x^{k,T^k - 1/2}$  satisfies condition~\eqref{hpe:monteiro-svaiter}, where $x^{k,T^k-1/2}$ is the output of \Cref{alg:eg}. This is exactly what we need. The following theorem provides an upper bound on the number of iterations $T^k$ required by \Cref{alg:eg} to terminate and produce the output $x_f^{k+1}$.
\begin{theorem}\label{eg:thm}
	Let $M$ satisfy
	\begin{equation}\label{eg:M}
		M \geq L_p.
	\end{equation}
	Then step~\eqref{eg:tensor_step} on \cref{eg:line:1} of \Cref{alg:eg} is well defined and the number of iterations $T^k$ performed by \Cref{alg:eg} is upper-bounded as follows:
	\begin{equation}\label{eg:T}
		T^k \leq \left(\lambda_k C_p(M,\sigma)\norm{x_g^k - x^{k,*}}^{p - 1} \right)^{2/p}+ 1,
	\end{equation}
	where $C_p$ is defined as
	\begin{equation}\label{eg:C}
		C_p(M,\sigma) = \frac{p^pM^p(1+\sigma^{-1})}{p!(pM - L_p)^{p/2}(pM+L_p)^{p/2-1}}.
	\end{equation}
\end{theorem}

\Cref{alg:eg} and \Cref{eg:thm} will further be used for the construction of the optimal high-order algorithm for solving problem~\eqref{eq:main}. It is worth mentioning the potential alternatives to \Cref{alg:eg} that we could use for gradient norm reduction.
For instance, we could use the tensor method of \citet{nesterov2021implementable}. However, the upper bound on the number of iterations for this method would involve the diameter of the level set of function $A_{\lambda_k}(\cdot;x_g^k)$ rather than the distance to the solution $\norm{x_g^{k} - x^{k,*}}$. This would be an obstacle towards development of the optimal algorithm. 
Alternatively, we could use the accelerated tensor method of \citet{nesterov2021implementable}. It turns out that it would work as we need. Moreover, the upper bound on the number of iterations would be even better than \eqref{eg:T}. However, we find the accelerated tensor method of \citet{nesterov2021implementable} to be too complicated, which could make the resulting optimal high-order method hard to implement. On the other hand, it would not give us any benefits for the construction of the optimal high-order method compared to \Cref{alg:eg}.

\subsection{Modification of the Analysis of A-HPE Framework}

Unfortunately, we cannot use \Cref{hpe:thm} for the analysis of our optimal algorithm. This is because inequality~\eqref{hpe:thm:eq} involves the distances $\norm{x_g^k - x_f^{k+1}}$ on the right-hand side. Hence, inequality~\eqref{hpe:thm:eq} does not allow us to estimate the iteration complexity $T^k$ of \Cref{alg:eg} using \Cref{eg:thm}. Further, we provide a new theorem that includes the analysis of the A-HPE framework and provides an upper bound on the distances $\norm{x_g^k - x^{k,*}}$.
\begin{theorem}\label{hpe:thm:2}
	The iterations of \Cref{alg:hpe} satisfy the following inequality:
	\begin{equation}
		2\beta_{K-1} (f(x_f^K) - f^*)  + \frac{1-\sigma}{1+\sigma}\sum_{k=0}^{K-1} \alpha_{k}^{-2}\sqn{x_g^k-x^{k,*}} \leq R^2.
	\end{equation}
\end{theorem}

\subsection{The First Optimal Tensor Method}

\begin{algorithm}[t]
	\caption{Optimal Tensor Method}
	\label{alg:ot}
	\begin{algorithmic}[1]
		\State {\bf input:} $x^0 = x_f^0 \in \R^d$
		\State {\bf parameters:} $\eta > 0, M > 0$, $\sigma \in (0,1)$, $K \in \{1,2,\ldots\}$
		\State $\beta_{-1} = 0$
		\For{$k=0,1,2,\ldots,K-1$}
		\State $\eta_k = \eta (1+k)^{(3p-1)/2}$
		\State $\beta_{k} = \beta_{k-1} + \eta_k$, $\lambda_k = \eta_k^2/\beta_k$, $\alpha_k = \eta_k / \beta_k$
		\State $x_g^k = \alpha_k x^k + (1-\alpha_k)x_f^k$
		\State $x^{k,0} = x_g^k$, $t=-1$
		\Repeat
		\State $t=t+1$
		\State $x^{k,t+1/2} = \argmin_{x \in \R^d} \Phi_{A_{\lambda_k}(\cdot;x_g^k)}^p (x;x^{k,t}) + \frac{pM}{(p+1)!}\norm{x - x^{k,t}}^{p+1}$
		\State $x^{k,t+1} = x^{k,t} - \left(\frac{M\norm{x^{k,t+1/2} - x^{k,t}}^{p-1}}{(p-1)!}\right)^{-1} \nabla A_{\lambda_k}(x^{k,t+1/2};x_g^k)$ 
		\Until{$\norm{\nabla A_{\lambda_k}(x^{k,t+1/2};x_g^k)} \leq \sigma\lambda_k^{-1}\norm{x^{k,t+1/2} - x^{k,0}}$}
		\State $T^k = t+1$
		\State $x_f^{k+1} = x^{k,T^k-1/2}$
		\State $x^{k+1} = x^k - \eta_k \nabla f(x_f^{k+1})$
		\EndFor
		\State {\bf output:} $x_f^K$
	\end{algorithmic}
\end{algorithm}

Now, we are ready to provide the first optimal high-order algorithm for solving problem~\eqref{eq:main}. In order to construct this algorithm, we use our Tensor Extragradient Method (\Cref{alg:eg}) to perform the computations on \cref{hpe:line:aux} of the A-HPE Framework (\Cref{alg:hpe}). We also use our choice of parameters $\eta_k$, $\beta_k$ and $\lambda_k$ which is provided by~\eqref{eq:eta} and~\eqref{eq:betalambda}. The resulting algorithm is formalized as \Cref{alg:ot}.

Now, we are ready to prove that \Cref{alg:ot} is an optimal algorithm. First, we need to establish an upper bound on the number of iterations $T^k$ performed by the inner repeat-loop of \Cref{alg:ot}. This is done by the following theorem.
\begin{theorem}\label{ot:thm}
	Let $M$ satisfy~\eqref{eg:M}. Then, the following inequality holds for \Cref{alg:ot}:
	\begin{equation}
		\sum_{k = 0}^{K-1} T^k \leq K +
		(1+K)
		\left(\frac{\eta(3p+1)^pC_p(M,\sigma)R^{p-1}}{2^p\sqrt{p}}\cdot 	\left(\frac{1+\sigma}{1-\sigma}\right)^{\frac{p-1}{2}}\right)^{\frac{2}{p}},
	\end{equation}
	where $C_p$ is defined by~\eqref{eg:C}.
\end{theorem}
\Cref{ot:thm} implies that with a proper choice of the parameter $\eta$, \Cref{alg:ot} performs $\cO(1)$ $p$-th order oracle calls per iteration on average. Indeed, let $\eta$ be chosen as follows:
\begin{equation}\label{ot:eta}
	\eta = \left(\frac{(3p+1)^pC_p(M,\sigma)R^{p-1}}{2^p\sqrt{p}}\cdot 	\left(\frac{1+\sigma}{1-\sigma}\right)^{\frac{p-1}{2}}\right)^{-1}.
\end{equation}
Then, \Cref{ot:thm} immediately implies
\begin{equation}\label{ot:T}
	\sum_{k = 0}^{K-1} T^k \leq 2K + 1.
\end{equation}
Finally, the following theorem establishes the total $p$-th order oracle complexity of \Cref{alg:ot}.
\begin{theorem}\label{ot:thm2}
	Let $M = L_p$ and $\sigma = 1/2$. Let $\eta$ be defined by~\eqref{ot:eta}. Then, to reach precision $f(x_f^k)  - f^* \leq \epsilon$, \Cref{alg:ot} requires no more than the following number of $p$-th order oracle calls:
	\begin{equation}
		5D_p\cdot\left(
		{L_pR^{p+1}}/{\epsilon}
		\right)^{\frac{2}{3p+1}} + 7,
	\end{equation}
	where $D_p$ is defined as follows:
	\begin{equation}
		D_p = 	\left(
		\frac{3^{\frac{p+1}{2}}(3p+1)^{p+1}p^p(p+1)}{2^{p+2}\sqrt{p}p!(p^2 - 1)^{\frac{p}{2}}}
		\right)^{\frac{2}{3p+1}}.
	\end{equation}
\end{theorem}

\Cref{ot:thm2} shows that the total $p$-th order oracle complexity of \Cref{alg:ot} is $\cO\left(\left(
{L_pR^{p+1}}/{\epsilon}
\right)^{\frac{2}{3p+1}} \right)$.
This oracle complexity matches the lower bounds of \citet{arjevani2019oracle} up to a universal constant that does not depend on $R$, $L_p$ and $\epsilon$. Hence, \Cref{alg:ot} is indeed the first optimal high-order algorithm for solving smooth convex minimization problems.

\bibliographystyle{apalike}
\bibliography{reference.bib}

\newpage

\appendix

\part*{Appendix}

\section{Proof of \Cref{lem:betalambda}}
The lower bound on $\beta_k$ can be obtained in the following way
\begin{align*}
	\beta_k &
	=
	\sum_{l=0}^{k}\eta (1 + l)^{\frac{3p-1}{2}}
	=
	\sum_{l=0}^{k}\eta \int_{0}^1 (1 + l)^{\frac{3p-1}{2}}dt
	\geq
	\sum_{l=0}^{k}\eta \int_{0}^1 (t + l)^{\frac{3p-1}{2}}dt
	\\&=
	\sum_{l=0}^{k}\eta \int_{l}^{l+1} t^{\frac{3p-1}{2}}dt
	=
	\eta \int_{0}^{k+1} t^{\frac{3p-1}{2}}dt
	=
	\frac{2\eta}{(3p+1)} (k+1)^{\frac{3p+1}{2}}.
\end{align*}
Upper bound on $\lambda_k$ is obtained using the lower bound on $\beta_k$ and~\eqref{hpe:aux3}.
\qed

\section{Proof of \Cref{eg:thm}}

\begin{align*}		
	\sqn{x^{k,t+1} - x^{k,*}}
	&=
	\sqn{x^{k,t} - x^{k,*}} + 2\<x^{k,t+1} - x^{k,t}, x^{k,t} - x^{k,*}> + \sqn{x^{k,t+1} - x^{k,t}}
	\\&=
	\sqn{x^{k,t} - x^{k,*}} + 2\<x^{k,t+1} - x^{k,t}, x^{k,t+1/2} - x^{k,*}>
	\\&
	+2\<x^{k,t+1} - x^{k,t}, x^{k,t} - x^{k,t+1/2}> + \sqn{x^{k,t+1} - x^{k,t}}
	\\&=
	\sqn{x^{k,t} - x^{k,*}} + 2\<x^{k,t+1} - x^{k,t}, x^{k,t+1/2} - x^{k,*}>
	\\&
	+\sqn{x^{k,t+1} - x^{k,t+1/2}} - \sqn{x^{k,t} - x^{k,t+1/2}}.
\end{align*}
From \eqref{eg:tensor_step} on\cref{eg:line:1} of \Cref{alg:eg}, we have
\begin{equation*}
	x^{k,t+1/2} = x^{k,t} -  \left(\frac{M\norm{x^{k,t+1/2} - x^{k,t}}^{p-1}}{(p-1)!}\right)^{-1} \nabla \Phi_{A^k}^p (x^{k,t+1/2};x^{k,t}).
\end{equation*}
Plugging this into the previous equation and using \cref{eg:line:2} of \Cref{alg:eg}, we get
\begin{align*}		
	\sqn{x^{k,t+1} - x^{k,*}}
	&=
	\sqn{x^{k,t} - x^{k,*}}
	- 2\gamma_{k,t}\<\nabla A^k(x^{k,t+1/2}), x^{k,t+1/2} - x^{k,*}>
	\\&
	+\gamma_{k,t}^2\sqn{\nabla \Phi_{A^k}^p (x^{k,t+1/2};x^{k,t}) - \nabla A^k(x^{k,t+1/2})} - \sqn{x^{k,t} - x^{k,t+1/2}},
\end{align*}
where $\gamma_{k,t} = \left(\frac{M\norm{x^{k,t+1/2} - x^{k,t}}^{p-1}}{(p-1)!}\right)^{-1} $.
Using the convexity of function $A^k(x)$, we get
\begin{align*}		
	\sqn{x^{k,t+1} - x^{k,*}}
	&=
	\sqn{x^{k,t} - x^{k,*}}
	- 2\gamma_{k,t}(A^k(x^{k,t+1/2}) - A^k(x^{k,*}))
	\\&
	+\gamma_{k,t}^2\sqn{\nabla \Phi_{A^k}^p (x^{k,t+1/2};x^{k,t}) - \nabla A^k(x^{k,t+1/2})} - \sqn{x^{k,t} - x^{k,t+1/2}}
	\\&\leq
	\sqn{x^{k,t} - x^{k,*}}
	+\gamma_{k,t}^2\sqn{\nabla \Phi_{A^k}^p (x^{k,t+1/2};x^{k,t}) - \nabla A^k(x^{k,t+1/2})}
	\\&- \sqn{x^{k,t} - x^{k,t+1/2}}.
\end{align*}
Using inequality~(1.6) of \citet{nesterov2021implementable}, we get
\begin{align*}		
	\sqn{x^{k,t+1} - x^{k,*}}
	&\leq
	\sqn{x^{k,t} - x^{k,*}}
	+\left(\frac{\gamma_{k,t}L_p}{p!}\right)^2\norm{x^{k,t+1/2} - x^{k,t}}^{2p} - \sqn{x^{k,t} - x^{k,t+1/2}}
	\\&=
	\sqn{x^{k,t} - x^{k,*}}
	- \left(1 - \left(\frac{\gamma_{k,t}L_p}{p!}\norm{x^{k,t+1/2} - x^{k,t}}^{p-1} \right)^2\right)\sqn{x^{k,t} - x^{k,t+1/2}}
	\\&=
	\sqn{x^{k,t} - x^{k,*}}
	- \left(1 - \left(\frac{L_p(p-1)!}{p!M} \right)^2\right)\sqn{x^{k,t} - x^{k,t+1/2}}
	\\&=
	\sqn{x^{k,t} - x^{k,*}}
	- \left(1 - \left(\frac{L_p}{pM} \right)^2\right)\sqn{x^{k,t} - x^{k,t+1/2}}.
\end{align*}
Using Lemma~1 of \citet{nesterov2021implementable}, we get
\begin{align*}		
	\sqn{x^{k,t+1} - x^{k,*}}
	&\leq
	\sqn{x^{k,t} - x^{k,*}}
	- \left(1 - \left(\frac{L_p}{pM} \right)^2\right)\left(\frac{p!}{pM + L_p}\norm{\nabla A^k(x^{k,t+1/2})}\right)^{2/p}.
\end{align*}
After telescoping and rearranging, for $T \leq T^k$ we get
\begin{align*}
	T\min_{t \in \{0,1,\ldots,T-1\}}\frac{(pM - L_p)(pM + L_p)}{p^2M^2}\left(\frac{p!}{pM + L_p}\norm{\nabla A^k(x^{k,t+1/2})}\right)^{2/p}
	\leq
	\sqn{x^{k,0} - x^{k,*}}.
\end{align*}
Taking both sides of the inequality in the power of $p/2$ gives
\begin{align*}
	\norm{x^{k,0} - x^{k,*}}^p
	&\geq
	T^{p/2}\min_{t \in \{0,1,\ldots,T-1\}}\frac{(pM - L_p)^{p/2}(pM + L_p)^{p/2}}{p^pM^p}\frac{p!}{pM + L_p}\norm{\nabla A^k(x^{k,t+1/2})}
	\\&=
	T^{p/2}\min_{t \in \{0,1,\ldots,T-1\}}\frac{p!(pM - L_p)^{p/2}(pM + L_p)^{p/2-1}}{p^pM^p}\norm{\nabla A^k(x^{k,t+1/2})}.
\end{align*}
After rearranging, we get
\begin{align*}
	\min_{t \in \{0,1,\ldots,T-1\}}\norm{\nabla A^k(x^{k,t+1/2})}
	\leq
	\frac{p^pM^p\norm{x^{k,0} - x^{k,*}}^p}{p!(pM - L_p)^{p/2}(pM + L_p)^{p/2-1}} \cdot \frac{1}{T^{p/2}}.
\end{align*}

Now, let us prove upper bound~\eqref{eg:T} by a contradiction. Suppose that~\eqref{eg:T} is not true. Hence,
\begin{equation*}
	T^k >\left( \frac{\lambda_kp^pM^p(1+\sigma^{-1})\norm{x^{k,0} - x^{k,*}}^{p - 1} }{p!(pM - L_p)^{p/2}(pM+L_p)^{p/2-1}}\right)^{2/p}+ 1.
\end{equation*}
This implies
\begin{align*}
	\min_{t \in \{0,1,\ldots,T-1\}}\left(\norm{\nabla A^k(x^{k,t+1/2})} - c\lambda_k^{-1}\norm{x^{k,0} - x^{k,*}} \right)\leq 0,
\end{align*}
where $c = (1+\sigma^{-1})^{-1}$ and $T = T^k - 1$.
Using the $\lambda_k^{-1}$-strong convexity of $A^k(x)$, we get
\begin{align*}
	0
	&\geq
	\min_{t \in \{0,1,\ldots,T-1\}}\left(\norm{\nabla A^k(x^{k,t+1/2})} - c\lambda_k^{-1}\norm{x^{k,0} - x^{k,*}} \right)
	\\&\geq
	\min_{t \in \{0,1,\ldots,T-1\}}\left(\norm{\nabla A^k(x^{k,t+1/2})}
	- c\lambda_k^{-1}\norm{x^{k,0} - x^{k,t+1/2}}
	-c\lambda_k^{-1}\norm{x^{k,*} - x^{k,t+1/2}}
	\right)
	\\&\geq
	\min_{t \in \{0,1,\ldots,T-1\}}\left(\norm{\nabla A^k(x^{k,t+1/2})}
	- c\lambda_k^{-1}\norm{x^{k,0} - x^{k,t+1/2}}
	-c\norm{\nabla A^k(x^{k,t+1/2})}
	\right)
	\\&=
	\min_{t \in \{0,1,\ldots,T-1\}}\left((1-c)\norm{\nabla A^k(x^{k,t+1/2})}
	- c\lambda_k^{-1}\norm{x^{k,0} - x^{k,t+1/2}}
	\right).
\end{align*}
Dividing by $1-c$ gives
\begin{align*}
	0
	&\geq
	\min_{t \in \{0,1,\ldots,T-1\}}\left(\norm{\nabla A^k(x^{k,t+1/2})}
	- \frac{\lambda_k^{-1}}{c^{-1}-1}\norm{x^{k,0} - x^{k,t+1/2}}
	\right).
\end{align*}
Plugging $c = (1+\sigma^{-1})^{-1}$ gives 
\begin{align*}
	0
	&\geq
	\min_{t \in \{0,1,\ldots,T-1\}}\left(\norm{\nabla A^k(x^{k,t+1/2})}
	- \sigma\lambda_k^{-1} \norm{x^{k,0} - x^{k,t+1/2}}
	\right).
\end{align*}
This means that the inner repeat-loop of \Cref{alg:eg} terminated after no more than $T$ iterations, which contradicts with $T^k = T+1$. This concludes the proof.
\qed

\section{Proof of \Cref{hpe:thm:2}}
Here we also provide the proof of \Cref{hpe:thm} for completeness.
Using \cref{hpe:line:grad} of \Cref{alg:hpe}, we get
\begin{align*}
	\sqn{x^{k+1} - x^*}
	&=
	\sqn{x^k - x^*} - 2\eta_k\<\nabla f(x_f^{k+1}),x^k - x^*> + \eta_k^2\sqn{\nabla f(x_f^{k+1})}.
\end{align*}
Using~\eqref{hpe:aux2}, we get $x^k = \alpha_k^{-1}x_g^k - \alpha_k^{-1}(1-\alpha_k)x_f^k$, which implies
\begin{align*}
	\sqn{x^{k+1} - x^*}
	&=
	\sqn{x^k - x^*} - 2\eta_k\<\nabla f(x_f^{k+1}),\alpha_k^{-1}x_g^k - \alpha_k^{-1}(1-\alpha_k)x_f^k- x^*> + \eta_k^2\sqn{\nabla f(x_f^{k+1})}
	\\&=
	\sqn{x^k - x^*} + \eta_k^2\sqn{\nabla f(x_f^{k+1})}
	\\&
	+2(\beta_k - \eta_k)\<\nabla f(x_f^{k+1}),x_f^k> 
	-2\beta_k\<\nabla f(x_f^{k+1}),x_g^k>
	+2\eta_k\<\nabla f(x_f^{k+1}),x^*>
	\\&=
	\sqn{x^k - x^*}
	+2(\beta_k - \eta_k)\<\nabla f(x_f^{k+1}),x_f^k-x_f^{k+1}> 
	+2\eta_k\<\nabla f(x_f^{k+1}),x^* - x_f^{k+1}>
	\\&
	-2\beta_k\<\nabla f(x_f^{k+1}),x_g^k - x_f^{k+1}>
	+ \eta_k^2\sqn{\nabla f(x_f^{k+1})}.
\end{align*}
Using the convexity of $f(x)$ and~\eqref{hpe:aux3}, we get
\begin{align*}
	\sqn{x^{k+1} - x^*}
	&\leq
	\sqn{x^k - x^*}
	+2(\beta_k - \eta_k)(f(x_f^k) - f(x_f^{k+1}))
	+2\eta_k( f^* - f(x_f^{k+1}))
	\\&
	-2\beta_k\<\nabla f(x_f^{k+1}),x_g^k - x_f^{k+1}>
	+ \eta_k^2\sqn{\nabla f(x_f^{k+1})}
	\\&=
	\sqn{x^k - x^*}
	-\beta_k (f(x_f^{k+1}) - f^*)
	+\beta_{k-1}(f(x_f^k) - f^*)
	\\&
	-2\beta_k\<\nabla f(x_f^{k+1}),x_g^k - x_f^{k+1}>
	+ \eta_k^2\sqn{\nabla f(x_f^{k+1})}.
\end{align*}
Using~\eqref{hpe:aux2}, we get
\begin{align*}
	\sqn{x^{k+1} - x^*}
	&\leq
	\sqn{x^k - x^*}
	-\beta_k (f(x_f^{k+1}) - f^*)
	+\beta_{k-1}(f(x_f^k) - f^*)
	\\&
	-2\<\eta_k\nabla f(x_f^{k+1}),\beta_k\eta_k^{-1}(x_g^k - x_f^{k+1})>
	+ \eta_k^2\sqn{\nabla f(x_f^{k+1})}
	\\&=
	\sqn{x^k - x^*}
	-\beta_k (f(x_f^{k+1}) - f^*)
	+\beta_{k-1}(f(x_f^k) - f^*)
	\\&
	+2\<\eta_k\nabla f(x_f^{k+1}),\alpha_k^{-1}(x_f^{k+1}-x_g^k)>
	+ \eta_k^2\sqn{\nabla f(x_f^{k+1})}.
\end{align*}
Using the parallelogram rule, we get
\begin{align*}
	\sqn{x^{k+1} - x^*}
	&\leq
	\sqn{x^k - x^*}
	-\beta_k (f(x_f^{k+1}) - f^*)
	+\beta_{k-1}(f(x_f^k) - f^*)
	\\&
	+\sqn{\eta_k\nabla f(x_f^{k+1}) + \alpha_k^{-1}(x_f^{k+1}-x_g^k)}
	-\alpha_k^{-2}\sqn{x_f^{k+1} - x_g^k}
	\\&=
	\sqn{x^k - x^*}
	-\beta_k (f(x_f^{k+1}) - f^*)
	+\beta_{k-1}(f(x_f^k) - f^*)
	\\&
	+\eta_k^2\sqn{\nabla f(x_f^{k+1}) + \eta_k^{-1}\alpha_k^{-1}(x_f^{k+1}-x_g^k)}
	-\alpha_k^{-2}\sqn{x_f^{k+1} - x_g^k}.
\end{align*}
Using~\eqref{hpe:aux2} and~\eqref{hpe:aux3}, we get
\begin{align*}
	\sqn{x^{k+1} - x^*}
	&\leq
	\sqn{x^k - x^*}
	-\beta_k (f(x_f^{k+1}) - f^*)
	+\beta_{k-1}(f(x_f^k) - f^*)
	\\&
	+\eta_k^2\sqn{\nabla f(x_f^{k+1}) + \lambda_k^{-1}(x_f^{k+1}-x_g^k)}
	-\alpha_k^{-2}\sqn{x_f^{k+1} - x_g^k}.
\end{align*}
Using~\eqref{hpe:aux}, we get
\begin{align*}
	\sqn{x^{k+1} - x^*}
	&\leq
	\sqn{x^k - x^*}
	-\beta_k (f(x_f^{k+1}) - f^*)
	+\beta_{k-1}(f(x_f^k) - f^*)
	\\&
	+\eta_k^2\lambda_k^{-2}\sigma^2\sqn{x_f^{k+1} - x_g^k}
	-\alpha_k^{-2}\sqn{x_f^{k+1} - x_g^k}
	\\&=
	\sqn{x^k - x^*}
	-\beta_k (f(x_f^{k+1}) - f^*)
	+\beta_{k-1}(f(x_f^k) - f^*)
	\\&
	-\alpha_k^{-2}(1-\sigma^2)\sqn{x_f^{k+1} - x_g^k}.
\end{align*} 
Now, let us bound $\norm{x_g^k - x^{k,*}}$ using $\lambda_k^{-1}$-strong convexity of $A_{\lambda_k} (\cdot;x_g^k)$ and~\eqref{hpe:aux}:
\begin{align*}
	\norm{x_g^k - x^{k,*}}
	&\leq
	\norm{x_g^k - x_f^{k+1}}
	+\norm{x_f^{k+1} - x^{k,*}}
	\\&\leq
	\norm{x_g^k - x_f^{k+1}}
	+\lambda_{k}\norm{\nabla A_{\lambda_k} (\cdot;x_g^k)}
	\\&\leq
	(1+\sigma)\norm{x_f^{k+1} - x_g^k}.
\end{align*}
Plugging this into the previous inequality gives
\begin{align*}
	\sqn{x^{k+1} - x^*}
	&\leq
	\sqn{x^k - x^*}
	-\beta_k (f(x_f^{k+1}) - f^*)
	+\beta_{k-1}(f(x_f^k) - f^*)
	\\&
	-\alpha_k^{-2}\frac{(1-\sigma^2)}{(1+\sigma)^2}\sqn{x_g^k - x^{k,*}}
	\\&=
	\sqn{x^k - x^*}
	-\beta_k (f(x_f^{k+1}) - f^*)
	+\beta_{k-1}(f(x_f^k) - f^*)
	\\&
	-\alpha_k^{-2}\frac{(1-\sigma)}{(1+\sigma)}\sqn{x_g^k - x^{k,*}}.
\end{align*} 
Rearranging and telescoping concludes the proof.
\qed

\section{Proof of \Cref{ot:thm}}
Using \Cref{eg:thm}, we get
\begin{align*}
	\left(\sum_{k=0}^{K-1} (T^k - 1)\right)^{\frac{p}{p-1}}
	&\leq		
	\left(\sum_{k=0}^{K-1} \left(\lambda_k C_p(M,\sigma)\norm{x_g^k - x^{k,*}}^{p - 1} \right)^{2/p} \right)^{\frac{p}{p-1}}.
\end{align*}
Let us choose parameters $\tau_0,\ldots,\tau_{K-1}$  as follows:
\begin{equation*}
	\tau_k = \left(\sum_{l=0}^{K-1}(1+l)^{p-1}\right)^{-1}(1+k)^{p-1}
\end{equation*}
Then, we have
\begin{align*}
	\left(\sum_{k=0}^{K-1} (T^k - 1)\right)^{\frac{p}{p-1}}
	&\leq		
	\left(\sum_{k=0}^{K-1} \tau_k\tau_k^{-1}\left(\lambda_k C_p(M,\sigma)\norm{x_g^k - x^{k,*}}^{p - 1} \right)^{2/p} \right)^{\frac{p}{p-1}}.
\end{align*}
Note that parameters $\tau_k$ satisfy
\begin{equation*}
	\sum_{k=0}^{K-1}\tau_k = 0, \qquad \tau_0,\ldots,\tau_{K-1} \geq 0.
\end{equation*}
Hence, using the convexity of function $(\cdot)^{p/(p-1)}$, we get
\begin{align*}
	\left(\sum_{k=0}^{K-1} (T^k - 1)\right)^{\frac{p}{p-1}}
	&\leq		
	\sum_{k=0}^{K-1}\tau_k\left( \tau_k^{-1}\left(\lambda_k C_p(M,\sigma)\norm{x_g^k - x^{k,*}}^{p - 1} \right)^{2/p} \right)^{\frac{p}{p-1}}
	\\&=
	C_p(M,\sigma)^{\frac{2}{p-1}}\sum_{k=0}^{K-1}\tau_k^{\frac{-1}{p-1}}\left(\lambda_k \right)^{\frac{2}{p-1}} \sqn{x_g^k - x^{k,*}}.
\end{align*}
Using \Cref{lem:betalambda}, we get
\begin{align*}
	\left(\sum_{k=0}^{K-1} (T^k - 1)\right)^{\frac{p}{p-1}}
	&\leq		
	C_p(M,\sigma)^{\frac{2}{p-1}}\sum_{k=0}^{K-1}\tau_k^{\frac{-1}{p-1}}\left(\frac{\eta(3p+1)}{2}(1+k)^{\frac{3(p-1)}{2}} \right)^{\frac{2}{p-1}} \sqn{x_g^k - x^{k,*}}
	\\&=
	\left(\frac{\eta(3p+1)C_p(M,\sigma)}{2} \right)^{\frac{2}{p-1}}\sum_{k=0}^{K-1}\tau_k^{\frac{-1}{p-1}}(1+k)^3 \sqn{x_g^k - x^{k,*}}.
\end{align*}
Using the definition of $\tau_k$, we get
\begin{align*}
	\left(\sum_{k=0}^{K-1} (T^k - 1)\right)^{\frac{p}{p-1}}
	&\leq		
	\left(\frac{\eta(3p+1)C_p(M,\sigma)}{2} \right)^{\frac{2}{p-1}}
	\cdot
	\left(\sum_{l=0}^{K-1}(1+l)^{p-1}\right)^{\frac{1}{p-1}}
	\\&\cdot
	\sum_{k=0}^{K-1}(1+k)^2 \sqn{x_g^k - x^{k,*}} .
\end{align*}
Using the inequality
\begin{align*}
	\sum_{l=0}^{K-1}(1+l)^{p-1}
	&\leq 
	\sum_{l=0}^{K-1} \int_{0}^1 (1+l + t)^{p-1} dt
	=
	\sum_{l=0}^{K-1} \int_{l+1}^{l+2} t^{p-1} dt
	\\&=
	\int_{1}^{K+1}  t^{p-1}dt = \frac{(1+K)^p - 1}{p}
	\leq \frac{1}{p} (1+K)^p,
\end{align*}
we get
\begin{align*}
	\left(\sum_{k=0}^{K-1} (T^k - 1)\right)^{\frac{p}{p-1}}
	&\leq		
	\left(\frac{\eta(3p+1)C_p(M,\sigma)}{2} \right)^{\frac{2}{p-1}}
	\cdot
	\left(\frac{1}{p} (1+K)^p\right)^{\frac{1}{p-1}}
	\\&\cdot
	\sum_{k=0}^{K-1}(1+k)^2 \sqn{x_g^k - x^{k,*}} 
	\\&=
	\left(\frac{\eta(3p+1)C_p(M,\sigma)}{2\sqrt{p}} \right)^{\frac{2}{p-1}}
	(1+K)^{\frac{p}{p-1}}
	\sum_{k=0}^{K-1}(1+k)^2 \sqn{x_g^k - x^{k,*}} .
\end{align*}
From \eqref{hpe:aux2} and \Cref{lem:betalambda}, we get
\begin{align*}
	\alpha_k^{-1} &= \frac{\beta_k}{\eta_k}
	\geq
	\frac{2\eta(1+k)^{\frac{3p+1}{2}}}{(3p+1)} \cdot\frac{1}{\eta(1+k)^{\frac{3p-1}{2}}}
	=
	\frac{2}{(3p+1)}(1+k).
\end{align*}
Hence,
\begin{align*}
	\left(\sum_{k=0}^{K-1} (T^k - 1)\right)^{\frac{p}{p-1}}
	&\leq
	\left(\frac{\eta(3p+1)C_p(M,\sigma)}{2\sqrt{p}} \right)^{\frac{2}{p-1}}
	(1+K)^{\frac{p}{p-1}}
	\sum_{k=0}^{K-1} \frac{(3p+1)^2}{4}\alpha_k^{-2}\sqn{x_g^k - x^{k,*}}
	\\&=
	\left(\frac{\eta(3p+1)C_p(M,\sigma)}{2\sqrt{p}} \right)^{\frac{2}{p-1}}
	\frac{(3p+1)^2(1+K)^{\frac{p}{p-1}}}{4}
	\sum_{k=0}^{K-1}\alpha_k^{-2}\sqn{x_g^k - x^{k,*}}.
\end{align*}
Using \Cref{hpe:thm:2}, we get
\begin{align*}
	\left(\sum_{k=0}^{K-1} (T^k - 1)\right)^{\frac{p}{p-1}}
	&\leq
	\left(\frac{\eta(3p+1)C_p(M,\sigma)}{2\sqrt{p}} \right)^{\frac{2}{p-1}}
	\frac{(3p+1)^2(1+K)^{\frac{p}{p-1}}}{4}
	\frac{1+\sigma}{1-\sigma}R^2.
\end{align*}
After taking both sides of the inequality in the power of $\frac{p-1}{p}$ we get
\begin{align*}
	\sum_{k=0}^{K-1} (T^k - 1)
	&\leq
	\left(\frac{\eta(3p+1)C_p(M,\sigma)}{2\sqrt{p}} \right)^{\frac{2}{p}}
	\frac{(3p+1)^{\frac{2(p-1)}{p}}(1+K)}{2^{\frac{2(p-1)}{p}}}
	\left(\frac{1+\sigma}{1-\sigma}R^2\right)^{\frac{p-1}{p}}
	\\&=
	(1+K)
	\left(\frac{\eta(3p+1)^pC_p(M,\sigma)R^{p-1}}{2^p\sqrt{p}}\cdot 	\left(\frac{1+\sigma}{1-\sigma}\right)^{\frac{p-1}{2}}\right)^{\frac{2}{p}}.
\end{align*}
After rearranging, we get
\begin{align*}
	\sum_{k=0}^{K-1} T^k
	&\leq
	K +
	(1+K)
	\left(\frac{\eta(3p+1)^pC_p(M,\sigma)R^{p-1}}{2^p\sqrt{p}}\cdot 	\left(\frac{1+\sigma}{1-\sigma}\right)^{\frac{p-1}{2}}\right)^{\frac{2}{p}}.
\end{align*}
\qed

\section{Proof of \Cref{ot:thm2}}

\Cref{hpe:thm:2} implies
\begin{align*}
	f(x_f^K) - f^*
	&\leq 
	R^2 /(2\beta_{K-1}).
\end{align*}
Using \Cref{lem:betalambda}, we get
\begin{align*}
	f(x_f^K) - f^*
	&\leq 
	\frac{(3p+1)R^2}{4\eta} \cdot\frac{1}{K^{\frac{3p+1}{2}}}.
\end{align*}
Choosing $K =\left\lceil\left(\frac{(3p+1)R^2}{4\eta\epsilon}\right)^{\frac{2}{3p+1}}\right\rceil$ implies  $f(x_f^K) - f^* \leq \epsilon$. Hence, we have the following upper bound on the total iteration complexity of \Cref{alg:ot}:
\begin{align*}
	K
	&\leq
	\left\lceil\left(\frac{(3p+1)R^2}{4\eta\epsilon}\right)^{\frac{2}{3p+1}}\right\rceil
	\\&\leq		
	\left(\frac{(3p+1)R^2}{4\eta\epsilon}\right)^{\frac{2}{3p+1}} + 1.
\end{align*}
Plugging $\eta$ defined by~\eqref{ot:eta} gives
\begin{align*}
	K
	&\leq
	\left(\frac{(3p+1)R^2}{4\epsilon}\cdot\frac{(3p+1)^pC_p(M,\sigma)R^{p-1}}{2^p\sqrt{p}}\cdot 	\left(\frac{1+\sigma}{1-\sigma}\right)^{\frac{p-1}{2}}\right)^{\frac{2}{3p+1}} + 1.
	\\&=
	\left(\frac{(3p+1)^{p+1}C_p(M,\sigma)R^{p+1}}{2^{p+2}\sqrt{p}\epsilon}\cdot 	\left(\frac{1+\sigma}{1-\sigma}\right)^{\frac{p-1}{2}}\right)^{\frac{2}{3p+1}} + 1.
\end{align*}
Using the definition of $C_p(M,\sigma)$, we get
\begin{align*}
	K
	&\leq
	\left(\frac{(3p+1)^{p+1}R^{p+1}}{2^{p+2}\sqrt{p}\epsilon}\cdot 	\left(\frac{1+\sigma}{1-\sigma}\right)^{\frac{p-1}{2}}\cdot\frac{p^pM^p(1+\sigma^{-1})}{p!(pM - L_p)^{p/2}(pM+L_p)^{p/2-1}}\right)^{\frac{2}{3p+1}}
	\\&+1.
\end{align*}
Using $M = L_p$, we get
\begin{align*}
	K
	&\leq
	\left(\frac{(3p+1)^{p+1}R^{p+1}}{2^{p+2}\sqrt{p}\epsilon}\cdot 	\left(\frac{1+\sigma}{1-\sigma}\right)^{\frac{p-1}{2}}\cdot\frac{p^pL_p(1+\sigma^{-1})}{p!(p-1)^{p/2}(p+1)^{p/2-1}}\right)^{\frac{2}{3p+1}}
	+1
	\\&=
	\left(
	\frac{L_pR^{p+1}}{\epsilon}\cdot
	\frac{(3p+1)^{p+1}p^p(p+1)}{2^{p+2}\sqrt{p}p!(p^2 - 1)^{\frac{p}{2}}}
	\cdot 
	\frac{(1+\sigma)^{\frac{p+1}{2}}}{\sigma(1-\sigma)^{\frac{p-1}{2}}}
	\right)^{\frac{2}{3p+1}}
	+1
	\\&=
	\left(
	\frac{L_pR^{p+1}}{\epsilon}
	\right)^{\frac{2}{3p+1}}
	\cdot
	\left(
	\frac{(3p+1)^{p+1}p^p(p+1)}{2^{p+2}\sqrt{p}p!(p^2 - 1)^{\frac{p}{2}}}
	\cdot 
	\frac{(1+\sigma)^{\frac{p+1}{2}}}{\sigma(1-\sigma)^{\frac{p-1}{2}}}
	\right)^{\frac{2}{3p+1}}
	+1.
\end{align*}
Using $\sigma = 1/2$, we get
\begin{align*}
	K
	&\leq
	\left(
	\frac{L_pR^{p+1}}{\epsilon}
	\right)^{\frac{2}{3p+1}}
	\cdot
	\left(
	\frac{(3p+1)^{p+1}p^p(p+1)}{2^{p+2}\sqrt{p}p!(p^2 - 1)^{\frac{p}{2}}}
	\cdot 
	3^{\frac{p+1}{2}}
	\right)^{\frac{2}{3p+1}}
	+1
	\\&=
	D_p\cdot\left(
	\frac{L_pR^{p+1}}{\epsilon}
	\right)^{\frac{2}{3p+1}} + 1.
\end{align*}
Finally, we have that \Cref{alg:ot} performs $(1 + 2T^k)$ of $p$-th order oracle calls at each iteration. Hence, using~\eqref{ot:T}, we get  following upper bound on the total oracle complexity:
\begin{align*}
	\sum_{k=0}^{K-1}(1+2T^k)
	&\leq K + 2(2K+1)
	= 5K + 2
	\\&\leq
	5D_p\cdot\left(
	\frac{L_pR^{p+1}}{\epsilon}
	\right)^{\frac{2}{3p+1}} + 7.
\end{align*}
\qed

\end{document}